\documentclass[a4paper]{scrartcl}

\newcommand{\hW}{\widehat W}
\usepackage[utf8]{inputenc}				%
\usepackage{amsmath}					%
\usepackage{mathtools}					%
\usepackage{amssymb,wasysym,amsthm}		%
\usepackage{amsfonts}					%
\usepackage{comment}					%
\usepackage[babel]{csquotes}			%
\usepackage{colonequals}				%
\usepackage{url}						%
\usepackage{graphicx}					%
\usepackage{bbm}						%
\usepackage{nicefrac}					%
\makeatletter
\let\@fnsymbol\@arabic
\makeatother
\theoremstyle{plain}
\newcounter{theoremCounter}
\numberwithin{theoremCounter}{section}

\newtheorem{proposition}[theoremCounter]{Proposition}
\newtheorem{theorem}[theoremCounter]{Theorem}

\theoremstyle{definition}

\newtheorem{remark}[theoremCounter]{Remark}

\newcommand{\id}{{\boldsymbol{\mathbbm{1}}}}
\newcommand{\eps}{\varepsilon}

\newcommand{\R}{\mathbb{R}}

\DeclareMathOperator{\PSym}{PSym}
\newcommand{\Rp}{\mathbb{R}_+}
\newcommand{\Rmat}[2]{\mathbb{R}^{#1\times #2}}
\newcommand{\Rnn}{\Rmat{n}{n}}

\newcommand{\PSymn}{\PSym(n)}

\DeclareMathOperator{\dev}{dev}
\DeclareMathOperator{\tr}{tr}

\newcommand{\intd}[1]{{\mathrm{d}#1}}

\newcommand{\dx}{\intd{x}}

\newcommand{\grad}{\nabla}

\newcommand{\norm}[1]{\Vert #1 \Vert}

\newcommand{\innerproduct}[1]{\langle #1 \rangle}

\newcommand{\iprod}{\innerproduct}

\renewcommand{\PSym}{\mathop{\mathrm{Sym_+}}}
\DeclareMathOperator{\D}{D}
\DeclareMathOperator{\Div}{Div}

\begin{document}
\title{\vspace*{-1.4em}On the convexity of nonlinear elastic energies in the right Cauchy-Green tensor}
\author{
David Yang Gao\thanks{School of Applied Science, Federation University Australia, Mt Helen VIC 3350, Australia, email: d.gao@federation.edu.au},\quad
Patrizio Neff\thanks{Patrizio Neff,  \ \ Head of Lehrstuhl f\"{u}r Nichtlineare Analysis und Modellierung, Fakult\"{a}t f\"{u}r Mathematik, Universit\"{a}t Duisburg-Essen,  Thea-Leymann Str. 9, 45127 Essen, Germany, email: patrizio.neff@uni-due.de},\quad
Ionel Roventa\thanks{Department of Mathematics, University of Craiova, A. I. Cuza Street, Nr. 13, 200585 Craiova, Romania, email: ionelroventa@yahoo.com},\quad
Christian Thiel\thanks{Corresponding author: Christian Thiel,  \ \ Lehrstuhl f\"{u}r Nichtlineare Analysis und Modellierung, Fakult\"{a}t f\"{u}r Mathematik, Universit\"{a}t Duisburg-Essen,  Thea-Leymann Str. 9, 45127 Essen, Germany, email: christian.thiel@uni-due.de } \quad
}
\date{\today\vspace*{-1.05em}}
\maketitle
\begin{abstract}
We present a sufficient condition under which a weak solution of the Euler-Lagrange equations in nonlinear elasticity is already a global minimizer of the corresponding elastic energy functional.
This criterion is applicable to energies $W(F)=\widehat{W}(F^TF)=\widehat{W}(C)$ which are convex with respect to the right Cauchy-Green tensor $C=F^TF$, where $F$ denotes the gradient of deformation.
Examples of such energies exhibiting a blow up for $\det F\to0$ are given.
\end{abstract}
\bigskip


\bigskip

\section{Introduction}
In a recent contribution, we have stated conditions for an elastic energy $W(F) = \hW(C)$ to be convex in terms of the right Cauchy-Green deformation tensor $C = F^TF$. These conditions include dependencies on the logarithm of the determinant, showing the necessary singular behavior for $\det F \to 0$; see \cite{lehmich2013convexity} for more details.

This subject has been taken up by Spector \cite{spector2015} and {\v{S}}ilhav{\'y} \cite{silhavy2015convexity}, who presented alternative proofs. The main result is the following. Here and throughout, $\PSymn$ denotes the set of positive definite symmetric $n\times n$--matrices and $\Rp\colonequals (0,\infty)$.
\begin{theorem}
\label{theorem:detConvexity}
Let $f\in C^2(\Rp)$. The function $C \mapsto f(\det C)$ is convex in $C$ on the set $\PSymn$ if and only if
\begin{equation}
	f''(s) + \frac{n-1}{n\,s}\,f'(s) \geq 0 \quad\text{ and }\quad f'(s)\leq0 \quad\text{ for all }\; s\in\Rp\,.
\end{equation}
In particular, these conditions are satisfied for $f(s)=-\log s$.
\end{theorem}
In this short note we would like to substantiate our claim from \cite{lehmich2013convexity} that convexity in $C$ should be somehow \enquote{nice}. It is clear that convexity of the strain energy in $C$ is independent of both the \emph{rank-one convexity} condition \cite{schroder2010poly} and Ball's \emph{polyconvexity} condition \cite{ball1977constitutive}. Therefore, convexity of $\hW$ as function of $C$ is not sufficient for establishing the existence of minimizers via the direct methods of the calculus of variations. Of course, in a neighborhood of the identity, an existence proof based on the implicit function theorem is always possible \cite{ciarlet2009pure}.

However, in \cite{lehmich2013convexity} we have given an example of a strain energy $W(F)=\hW(C)$ which is both polyconvex with respect to $F$ and convex in $C$; furthermore, the energy shows the correct behaviour for infinite compression, i.e.\ $\hW(C) \to \infty$ as $\det C \to 0$. Convexity with respect to $C$ is also used profitably in the derivation of some models in elasto-plasticity, see \cite{gupta2007evolution}.
\section{Applications of the convexity with respect to $C$}
We consider the problem of minimizing the elastic energy $I(\varphi)=\int_\Omega W(\grad\varphi)\,\intd x$, where $\varphi\colon \Omega \to \R^3$ is the deformation of an elastic body $\Omega\subset\R^3$, subject to the boundary condition $\varphi|_\Gamma = \phi_0$ with $\Gamma\subset\partial\Omega$, see \cite{Ciarlet1988}. The corresponding formal Euler-Lagrange equation is $0 = \Div \D_F\! W(F) = \Div S_1(F)$, where $F = \nabla \varphi$ denotes the gradient of deformation, $S_1(F) = \D_F\! W(F) = F\,S_2(F^TF)$ is the first Piola-Kirchhoff stress tensor and $S_2(C)=2\,\D_C\! \hW(C)$ is the second Piola-Kirchhoff stress tensor.
The following simple observation extracted from \cite{gao1992global,gao2008closed} shows the applicability of our convexity condition.

\begin{proposition}\label{prop:1}
Assume that $\hW$ is convex in $C$ and that $\varphi_0$ is a weak solution of the corresponding Euler-Lagrange equation.
If the second Piola-Kirchhoff stress tensor $S_2(C_0) = 2\,\D_C\! \hW(C_0)$ corresponding to $C_0=\nabla \varphi_0^T\nabla \varphi_0$ is positive definite everywhere in $\Omega$, then this solution is a global minimizer.
\end{proposition}\vspace*{.5em}
\begin{remark}
Note that the condition of positive definiteness only refers to the second Piola-Kirchhoff stress \emph{evaluated at the weak solution $\varphi_0$} of the Euler-Lagrange equation.
\end{remark}\vspace*{.5em}
Convexity with respect to $C$ is tantamount to the monotonicity of the second Piola-Kirchhoff stress tensor $S_2$, i.e.\ to the condition
\[
	\iprod{S_2(C_1)-S_2(C_2),\, C_1-C_2} \;\geq\; 0\,,
\]
where $\iprod{X,Y}=\tr(Y^TX)$ denotes the canonical inner product on $\Rnn$. This does not violate any physical principle, whereas monotonicity of the first Piola-Kirchhoff stress tensor $S_1$ (i.e.\ convexity with respect to the deformation gradient $F$) must be excluded to ensure physically plausible material behaviour.

The positive definiteness of the second Piola-Kirchhoff stress $S_2(C_0)$ at $\varphi_0$, on the other hand, means that all occurring forces, transformed to the reference configuration, are of a \enquote{tensional} type, i.e.\ that the solution $\varphi_0$ must not exhibit any \enquote{compressional} forces. Note, again, that this condition only refers to the solution $\varphi_0$ and that we do not require $S_2(C)=2\,\D_C\! \hW(C)$ to be convex for \emph{all} $C\in\PSym(3)$, which together with the convexity of $\hW$ in $C$ would imply convexity of $W$ in $F$.
\begin{proof}
We write $\varphi$ as $\varphi = \varphi_0 + u$. Then $\nabla \varphi = \nabla \varphi_0 + \nabla u$ and thus
\[
	\nabla \varphi^T\nabla \varphi\ =\ (\nabla \varphi_0 + \nabla u)^T(\nabla \varphi_0 + \nabla u)\ =\ \nabla \varphi_0^T\nabla \varphi_0 + \nabla \varphi_0^T\nabla u + \nabla u^T\nabla \varphi_0 + \nabla u^T\nabla u\,.
\]
We therefore find $\nabla \varphi^T\nabla \varphi - \nabla \varphi_0^T\nabla \varphi_0\ =\ \nabla \varphi_0^T\nabla u + \nabla u^T\nabla \varphi_0 + \nabla u^T\nabla u$.
The convexity of $\hW$ with respect to $C$ implies
\begin{equation}\label{eq:1.1}
\hW(C) - \hW(C_0)\ \geq\ \innerproduct{\D_C\!\hW(C_0), C - C_0} \qquad\text{for all }\; C \in \PSymn\,.
\end{equation}
Using \eqref{eq:1.1} we obtain the following minimizing estimates:
\begin{align}
&\int_\Omega \hW(\nabla \varphi^T\nabla \varphi) - \hW(\nabla \varphi_0^T\nabla \varphi_0)\,\intd x\nonumber\\
&\qquad\qquad\geq\ \int_\Omega \innerproduct{\D_C\!\hW(\nabla \varphi_0^T\nabla \varphi_0), \nabla \varphi_0^T\nabla u + \nabla u^T\nabla \varphi_0 + \nabla u^T\nabla u}\,\intd x\nonumber\\
&\qquad\qquad=\ \int_\Omega 2\innerproduct{S_1(\nabla \varphi_0), \nabla u} + \frac12\,\innerproduct{S_2(C_0), \nabla u^T\nabla u}\,\intd x\nonumber\\
&\qquad\qquad=\ \frac12\,\int_\Omega \innerproduct{S_2(C_0), \nabla u^T\nabla u}\,\intd x \ \geq\ 0\,, \label{eq:eulerLagrangeComputationLastLine}
\end{align}
where the equality and inequality in \eqref{eq:eulerLagrangeComputationLastLine} follow from the fact that $\int_\Omega 2\innerproduct{S_1(\nabla \varphi_0), \nabla u}\dx = 0$ due to $\varphi_0$ satisfying the weak form of Euler-Lagrange equation and the assumed positive definiteness of $S_2(C_0)$, respectively.
\end{proof}
If the convexity condition in Proposition \ref{prop:1} is strengthened to the uniform convexity requirement $D_C^2 \hW(C).(H,H)\geq c^+\,\norm{H}^2$ for all symmetric matrices $H$ and some $c^+>0$, and if the solution $\varphi_0$ is a diffeomorphism, then it is easy to show that
\[
	\int_\Omega D_F^2 W(\grad\varphi_0).(\grad\vartheta,\grad\vartheta)\,\dx \;\geq\; \widetilde{c}^+ \,\norm{\grad\vartheta}^2_{L^2(\Omega)}
\]
for all test functions $\vartheta\in C_0^\infty(\Omega)$ and some $\widetilde{c}^+>0$, implying the local stability of the solution \cite{neff2002korn}.

We will apply Proposition \ref{prop:1} to a modified version of $W$ as given in \cite{lehmich2013convexity}.
\begin{proposition}\label{prop:2}
Let
\[
	W(F) = \hW(C) = \hW(F^TF) = \alpha\bigl(\tr(C)\bigr)^2 + \beta \tr(C^2) - \log \det C
\]
with $6\alpha + 2\beta - 1 = 0$ (which is strictly convex in $C$ and polyconvex provided that $\alpha, \beta > 0$). Assume that $\hW(C)$ admits a smooth solution $\varphi_0$ of the Euler-Lagrange equation satisfying $\det \nabla \varphi_0 \geq C^+ > 0$ and $S_2(\nabla \varphi_0)$ is positive definite everywhere. Then $\varphi_0$ is a global minimizer.
\end{proposition}

The condition $6\alpha + 2\beta = 1$, $\alpha, \beta > 0$ ensures that the condition $\D_C\hW(\id) = 0$ for a stress free reference configuration is fulfilled.
Similarly, for the energy $\hW(C) = \alpha \tr(C) + \beta \tr(C^2) - \log \det C$, we need to impose the conditions\; $\alpha + 2\,\beta = 1$,\; $\alpha>0$\; and\; $\beta>0$.
Since the energy $W$ given in Proposition \ref{prop:2} is polyconvex and coercive in the Sobolev space $W^{1,4}$ if $\alpha > 0$, it follows from the direct methods of the calculus of variations \cite{Ciarlet1988} that there exists a global minimizer to the boundary value problem considered in Proposition \ref{prop:1}.
The conditions under which a global minimizer satisfies the weak form of Euler-Lagrange equation are, however, not clear since it is not known whether $\det \nabla \varphi_0 \geq c^+$ for some $c^+ > 0$. On the other hand, most computational methods will try to solve the Euler-Lagrange equations. Thus, in general, it is not known whether a global minimizer has been found.

The energy
\[
	\hW(C) = \mu\Bigl(\alpha \bigl(\tr(C)\bigr)^2 + \beta\tr (C) - \log \det C\Bigr) %
\]
is a special member of the compressible Neo-Hooke family. It generates the linear elastic response $\mu\,\norm{\dev\eps}^2 + \mu\,\Big(2\,\alpha + \frac{1}{3}\Big)\,[\tr(\eps)]^2$ with an induced linear Poisson coefficient $\nu = \frac12\,\frac{4\,\alpha}{4\,\alpha + 1}\in(0,\frac12)$. Here, $\|\,.\,\|$ denotes the Frobenius matrix norm and $\eps$ is the linearized strain tensor.
Our examples show that the search for further insight into the properties of specific elastic energies might lead to unexpected new results. Certainly, the class of convex energies in $C$ is much to narrow from an application point of view. Since, however, all invariance principles and side conditions have been respected, our examples show that there may be enough space for more discoveries.

A more general application of the convexity in $C$ was given by Le Dret and Raoult \cite{ledret1995quasiconvexInBandle}: they showed that if $W(F)=\hW(C)$ is convex with respect to $C$, then the \emph{quasiconvex hull} $QW$ of $W$ is given by an analogue of \emph{Pipkin's formula} \cite{pipkin1994relaxed}:
\begin{equation}
	QW(F) = \inf_{S\in\PSym(3)} \hW(F^TF+S)\,.
\end{equation}
\section{Further examples}
A trivial example of a function $\hW$ which is convex in $C$ is the Saint Venant-Kirchhoff energy \cite{raoult1986non,ledret1995svk,gao2015analytic}
\[
\hW_{\mathrm{SVK}}(C)\ =\ \frac \mu 4\|C-\id\|^2 + \frac \lambda 8\tr(C - \id)^2\,,
\]
which is convex in $C$ as long as $\mu>0$, $3\lambda + 2\mu > 0$.
However, $\hW_{\mathrm{SVK}}$ is not rank-one convex, as has been shown by Raoult \cite{raoult1986non}.

For another class of examples, consider an energy function satisfying the \emph{Valanis-Landel hypothesis}, i.e.\ an energy of the form
\[
	W(F) = w(\lambda_1) + w(\lambda_2) + w(\lambda_3)
\]
where $w\colon \Rp\to\R$ is a fixed function and $\lambda_1,\lambda_2,\lambda_3$ are the singular values of $F$, i.e.\ the eigenvalues of $U=\sqrt{C}$. Since the eigenvalues of $U$ are the square roots of the eigenvalues of $C$, the function $W$ can be represented in terms of the eigenvalues of $C$ via $W(F) = \widetilde{w}(\lambda_1^2) + \widetilde{w}(\lambda_2^2) + \widetilde{w}(\lambda_3^2)$ with $\widetilde{w}(t)=w(\sqrt{t})$. %
Then according to a result by Davis \cite{davis1957}, $W$ is convex with respect to $C$ if and only if $\widetilde{w}$ is convex. If $f\colon \Rp\to\R$ satisfies the requirements of Theorem \ref{theorem:detConvexity}, then an energy of the form
\[
	W(F) = \widetilde{w}(\lambda_1^2) + \widetilde{w}(\lambda_2^2) + \widetilde{w}(\lambda_3^2) \;+\; g(\det C)
\]
is convex in $C$ as well if $\widetilde{w}$ is convex.
\section*{Acknowledgement}
The third author was supported by Grant PN-II-RU-TE-2014-4-1109 of the Romanian National Authority for Scientific Research, CNCS-UEFISCDI.

\begin{thebibliography}{10}

\bibitem{ball1977constitutive}
J.~Ball.
\newblock Constitutive inequalities and existence theorems in nonlinear
  elastostatics.
\newblock In {\em Nonlinear Analysis and Mechanics: Heriot-Watt Symposium},
  volume~1, pages 187--241. Pitman Publishing Ltd. Boston, 1977.

\bibitem{Ciarlet1988}
P.~G. Ciarlet.
\newblock {\em Three-Dimensional Elasticity}.
\newblock Number~1 in Studies in mathematics and its applications. Elsevier
  Science, 1988.

\bibitem{ciarlet2009pure}
P.~G. Ciarlet and C.~Mardare.
\newblock The pure displacement problem in nonlinear three-dimensional
  elasticity: intrinsic formulation and existence theorems.
\newblock {\em Comptes Rendus Mathematique}, 347(11):677--683, 2009.

\bibitem{davis1957}
C.~Davis.
\newblock All convex invariant functions of hermitian matrices.
\newblock {\em Archiv der Mathematik}, 8(4):276--278, 1957.

\bibitem{gao1992global}
D.~Y. Gao.
\newblock Global extremum criteria for finite elasticity.
\newblock {\em Zeitschrift f{\"u}r angewandte Mathematik und Physik},
  43(5):924--937, 1992.

\bibitem{gao2015analytic}
D.~Y. Gao and E.~Hajilarov.
\newblock Analytic solutions to three-dimensional finite deformation problems
  governed by {S}t {V}enant--{K}irchhoff material.
\newblock {\em Mathematics and Mechanics of Solids}, 2015.
\newblock DOI: 10.1177/1081286515591084.

\bibitem{gao2008closed}
D.~Y. Gao and R.~W. Ogden.
\newblock Closed-form solutions, extremality and nonsmoothness criteria in a
  large deformation elasticity problem.
\newblock {\em Zeitschrift f{\"u}r angewandte Mathematik und Physik},
  59(3):498--517, 2008.

\bibitem{gupta2007evolution}
A.~Gupta, D.~J. Steigmann, and J.~S. St{\"o}lken.
\newblock On the evolution of plasticity and incompatibility.
\newblock {\em Mathematics and Mechanics of Solids}, 12(6):583--610, 2007.

\bibitem{ledret1995svk}
H.~Le~Dret and A.~Raoult.
\newblock The quasiconvex envelope of the {S}aint {V}enant--{K}irchhoff stored
  energy function.
\newblock {\em Proceedings of the Royal Society of Edinburgh: Section A
  Mathematics}, 125(06):1179--1192, 1995.

\bibitem{ledret1995quasiconvexInBandle}
H.~Le~Dret and A.~Raoult.
\newblock Quasiconvex envelopes of stored energy densities that are convex with
  respect to the strain tensor.
\newblock In C.~Bandle, M.~Chipot, J.~S.~J. Paulin, J.~Bemelmans, and
  I.~Shafrir, editors, {\em Calculus of Variations, Applications and
  Computations}, volume 326, pages 138--146. CRC Press, 1995.

\bibitem{lehmich2013convexity}
S.~Lehmich, P.~Neff, and J.~Lankeit.
\newblock On the convexity of the function ${C} \mapsto f (\det {C})$ on
  positive-definite matrices.
\newblock {\em Mathematics and Mechanics of Solids}, 19(4):369--375, 2014.

\bibitem{neff2002korn}
P.~Neff.
\newblock On {K}orn's first inequality with non-constant coefficients.
\newblock {\em Proceedings of the Royal Society of Edinburgh: Section A
  Mathematics}, 132(01):221--243, 2002.

\bibitem{pipkin1994relaxed}
A.~C. Pipkin.
\newblock Relaxed energy densities for large deformations of membranes.
\newblock {\em IMA Journal of Applied Mathematics}, 52(3):297--308, 1994.

\bibitem{raoult1986non}
A.~Raoult.
\newblock Non-polyconvexity of the stored energy function of a {S}aint
  {V}enant--{K}irchhoff material.
\newblock {\em Aplikace matematiky}, 31(6):417--419, 1986.

\bibitem{schroder2010poly}
J.~Schr{\"o}der and P.~Neff.
\newblock {\em Poly-, quasi- and rank-one convexity in applied mechanics},
  volume 516.
\newblock Springer Science \& Business Media, 2010.

\bibitem{silhavy2015convexity}
M.~{\v{S}}ilhav{\'y}.
\newblock The convexity of {$C\mapsto h(\det C)$}.
\newblock {\em Technische Mechanik}, 35(1):60--61, 2015.

\bibitem{spector2015}
S.~J. Spector.
\newblock A note on the convexity of {$C\mapsto h(\det C)$}.
\newblock {\em Journal of Elasticity}, 118(2):251--256, 2015.

\end{thebibliography}
\end{document}